\documentclass[12pt,twoside]{article}

\begin{document}

\date{}





\centerline{}

\centerline {\Large{\bf On Fixed-point theorems in Intuitionistic  }}

\centerline{}

\centerline {\Large{\bf Fuzzy metric Space I}}

\centerline{}


\centerline{\bf {T.K. Samanta and Sumit Mohinta }}

\centerline{}

\centerline{Department of Mathematics, Uluberia College, India-711315.} %
\centerline{e-mail: mumpu$_{-}$tapas5@yahoo.co.in}
\centerline{e-mail: sumit.mohinta@yahoo.com}

\bigskip
\centerline{\bf Abstract}
 {\emph{In this paper, first we have
established two sets of sufficient conditions for a TS-IF
contractive mapping to have unique fixed point in a intuitionistic
fuzzy metric space. Then we have defined \,$(\,\varepsilon \,,\,
\lambda\,)$\, IF-uniformly locally contractive mapping and
\,$\eta\,-$\,chainable space, where it has been proved that the
\,$(\,\varepsilon \,,\, \lambda\,)$\, IF-uniformly locally
contractive mapping possesses a fixed point.}}

\centerline{}

\newtheorem{Theorem}{\quad Theorem}[section]

\newtheorem{Definition}[Theorem]{\quad Definition}

\newtheorem{Corollary}[Theorem]{\quad Corollary}

\newtheorem{Lemma}[Theorem]{\quad Lemma}

\newtheorem{Note}[Theorem]{\quad Note}

\newtheorem{Remark}[Theorem]{\quad Remark}

\newtheorem{Result}[Theorem]{\quad Result}

\newtheorem{Proposition}[Theorem]{\quad Proposition}

\newtheorem{Example}[Theorem]{\quad Example}

{\bf 2000 MSC No:} 03F55, 46S40.


\section{Introduction}

Fuzzy set theory was first introduce by Zadeh\cite{zadeh} in 1965 to describe
the situation in which data are imprecise or vague or uncertain.
Thereafter the concept of fuzzy set was generalized as intuitionistic
fuzzy set by K. Atanassov\cite{Atanasov,Atanasov1} in 1984. It has a wide
range of application in the field of population dynamics , chaos control ,
computer programming , medicine , etc.\\
The concept of fuzzy metric was first introduced by Kramosil and
Michalek\cite{Kramosil} but using the idea of intuitionistic fuzzy set,
Park\cite{Park} introduced the notion of intuitionistic fuzzy metric spaces
with the help of continuous t-norms and continuous t-conorms, which is a
generalization of fuzzy metric space due to George and Veeramani\cite{Veeramani}. \\
Introducing the contraction mapping with the help of the membership
function for fuzzy metric, several
authors\cite{Abdul,Grab,Veeramani2} established the Banach fixed
point theorem in fuzzy metric space. In the paper\cite{Abdul}, to
prove the Banach Fixed Point theorem in intuitionistic fuzzy metric
space, Mohamad\cite{Abdul} also introduced one concept of
contractive mapping, which is not so natural. There he proved that
every iterative sequence is a contractive sequence and then assumed
that every contractive sequences are Cauchy. But all these
contraction mappings, which they have considered to establish
different type fixed point theorem, do not bear the intension of the
contraction mapping with respect to a fuzzy metric, when a fuzzy
metric gives the degree of nearness of two points with respect to a
parameter $t$. Considering this meaning of fuzzy metric, in our
paper\cite{T.Samanta}, we have redefined the notion of contraction
mapping in a intuitionistic fuzzy metric space and then directly, it
has been proved that the every iterative sequence is a Cauchy
sequence, that is, we don't need to assume that every contractive
sequences are Cauchy sequences. Thereafter we have established the
Banach Fixed Point theorem there. In this paper, first we have
established two sets of sufficient conditions for a TS-IF
contractive mapping to have unique fixed point in a intuitionistic
fuzzy metric space. Then we have defined \,$(\,\varepsilon \,,\,
\lambda\,)$\, IF-uniformly locally contractive mapping and
\,$\eta\,-$\,chainable space, where it has been proved that  the
\,$(\,\varepsilon \,,\, \lambda\,)$\, IF-uniformly locally
contractive mapping possesses a fixed point.


\section{Preliminaries}

We quote some definitions and statements of a few theorems which will be
needed in the sequel.

\begin{Definition} \cite{Schweizer}.
A binary operation \, $\ast \; : \; [\,0 \; , \; 1\,] \; \times \;
[\,0 \; , \; 1\,] \;\, \longrightarrow \;\, [\,0 \; , \; 1\,]$ \, is
continuous \, $t$ - norm if \,$\ast$\, satisfies the
following conditions \, $:$ \\
$(\,i\,)$ \hspace{0.2cm} $\ast$ \, is commutative and associative ,
\\ $(\,ii\,)$ \hspace{0.1cm} $\ast$ \, is continuous , \\
$(\,iii\,)$ \hspace{0.01cm} $a \;\ast\;1 \;\,=\;\, a \hspace{1.2cm}
\forall \;\; a \;\; \varepsilon \;\; [\,0 \;,\; 1\,]$ , \\
$(\,iv\,)$ \hspace{0.1cm} $a \;\ast\; b \;\, \leq \;\, c \;\ast\; d$
\, whenever \, $a \;\leq\; c$  ,  $b \;\leq\; d$  and  $a \, , \, b
\, , \, c \, , \, d \;\, \varepsilon \;\;[\,0 \;,\; 1\,]$.
\end{Definition}

\begin{Definition}
\cite{Schweizer}. A binary operation \, $\diamond \; : \; [\,0 \; ,
\; 1\,] \; \times \; [\,0 \; , \; 1\,] \;\, \longrightarrow \;\,
[\,0 \; , \; 1\,]$ \, is continuous \, $t$-conorm if \,$\diamond$\,
satisfies the
following conditions \, $:$ \\
$(\,i\,)\;\;$ \hspace{0.01cm} $\diamond$ \, is commutative and
associative ,
\\ $(\,ii\,)\;$ \hspace{0.01cm} $\diamond$ \, is continuous , \\
$(\,iii\,)$ \hspace{0.01cm} $a \;\diamond\;0 \;\,=\;\, a
\hspace{1.2cm}
\forall \;\; a \;\; \in\;\; [\,0 \;,\; 1\,]$ , \\
$(\,iv\,)$ \hspace{0.1cm} $a \;\diamond\; b \;\, \leq \;\, c
\;\diamond\; d$ \, whenever \, $a \;\leq\; c$  ,  $b \;\leq\; d$
 and  $a \, , \, b \, , \, c \, , \, d \;\; \in\;\;[\,0
\;,\; 1\,]$.
\end{Definition}

\begin{Result}
\cite{klement}. $(\,a\,)\;$  For any \, $r_{\,1} \; , \;
r_{\,2} \;\; \in\;\; (\,0 \;,\; 1\,)$ \, with \, $r_{\,1} \;>\;
r_{\,2}$, there exist $\;r_{\,3} \; , \; r_{\,4} \;\; \in \;\; (\,0
\;,\; 1\,)$ \, such that \, $r_{\,1} \;\ast\; r_{\;3} \;>\; r_{\,2}$
\, and \, $r_{\,1} \;>\; r_{\,4} \;\diamond\; r_{\,2}.$
\\ $(\,b\,)$ \, For any \, $r_{\,5} \;\,
\in\;\, (\,0 \;,\; 1\,)$ , there exist \, $r_{\,6} \; , \; r_{\,7}
\;\, \in\;\, (\,0 \;,\; 1\,)$ \, such that \, $r_{\,6} \;\ast\;
r_{\,6} \;\geq\; r_{\,5}$ \,and\, $r_{\,7} \;\diamond\; r_{\,7}
\;\leq\; r_{\,5}.$
\end{Result}

\begin{Definition}
\cite{Park}
 Let \,$\ast$\, be a continuous \,$t$-norm ,
\,$\diamond$\, be a continuous \,$t$-conorm and $X$ be any non-empty set.
An \textbf{intuitionistic fuzzy metric} or in short $\textit{\textbf{IFM}}$ on
\,$X$\, is an object of the form \\ $A \;\,=\;\, \{\; (\,(\,x \;,\,y \;,\;
t\,) \;,\; \mu\,(\,x \;,\,y\;,\; t\,) \;,\; \nu\,(\,x \;,\,y \;,\; t\,) \;) \;\,
: \;\, (\,x \;,\,y \;,\; t\,) \;\,\in\;\, X^{2}\times(0\,,\,\infty)\}$\,
where $\mu \,,\, \nu\;$ are fuzzy sets on \, $X^{2}\times(0\,,\,\infty)$ ,
\,$\mu$\, denotes the degree of nearness and
\,$\nu$\, denotes the degree of non$-$nearness of $x$ and $y$ relative
to $t$ satisfying the following conditions $:$ \, for all $x , y , z
\, \in \,X , \, s , t \, > \, 0$ \\
$(\,i\,)$ \hspace{0.10cm}  $\mu\,(\,x \;,y\;, t\,) \;+\; \nu\,(\,x
\;,y\;, t\,) \;\,\leq\;\, 1 \hspace{0.5cm} \forall \;\; (\,x \;,\,y \;,\;
t\,)
\;\,\in\;\,X^{2}\times(0,\infty) ;$
\\$(\,ii\,)$ \hspace{0.10cm}$\mu\,(\,x \;,\,y \;,\; t\,) \;\,>\;\, 0 \, ;$ \\
$(\,iii\,)$ $\mu\,(\,x \;,\,y\;,\; t\,) \;\,=\;\, 1$ \, if
and only if \, $x \;=\;y \,$\,
\\$(\,iv\,)$ $\mu\,(\,x \;,\,y\;,\; t\,) \;\,=\;\,\mu\,(\,y \;,\,x\;,\; t\,);$
\\ $(\,v\,)$
\hspace{0.01cm} $\mu\,(\,x \;,\,y \;,\; s\,) \;\ast\; \mu\,(\,y \;,\,z \;,\; t\,)
\;\,\leq\;\, \mu\,(\,x \;,\,z \;,\; s\;+\;t\,
\,) \, ;$
\\$(\,vi\,)$ $\mu\,(\,x \,,\,y\,,\,\cdot\,) :(0 \,,\;\infty\,)\,\rightarrow
\,(0 \,,\;1]$ \, is continuous;
\\$(\,vii\,)$ \hspace{0.10cm}$\nu\,(\,x \;,\,y \;,\; t\,) \;\,>\;\, 0 \, ;$
\\$(\,viii\,)$ $\nu\,(\,x \;,\,y\;,\; t\,) \;\,=\;\, 0$ \, if
and only if \, $x \;=\;y \,;$\,
\\$(\,ix\,)$ $\nu\,(\,x \;,\,y\;,\; t\,) \;\,=\;\,\nu\,(\,y \;,\,x\;,\; t\,);$
\\ $(\,x\,)$\hspace{0.10cm} $\nu\,(\,x \;,\,y \;,\; s\,) \;\diamond\;
\mu\,(\,y \;,\,z \;,\; t\,)
\;\,\geq\;\, \nu\,(\,x \;,\,z \;,\; s\;+\;t\,
\,) \, ;$
\\$(\,xi\,)$ $\nu\,(\,x \,,\,y\,,\,\cdot\,) :(0\,,\;\infty\,)\,\rightarrow
\,(0\,,\;1]$ \, is continuous.
\\\\ If \,$A$\, is a $\textit{\textbf{IFM}}$ on
\,$X$ , the pair \,$(\,X \,,\, A\,)$\, will be called a \textbf{intuitionistic
fuzzy metric space} or in short $\textit{\textbf{IFMS}}$.\\\\ We further assume
that \,$(\,X \,,\, A\,)$\, is a $\textit{\textbf{IFMS}}$ with the property \\
$(\,xii\,)\;$ For all \,$a\,\in\,(0 \,,\, 1)$,\, $a \,\ast\, a \,=\, a$\,
and \,$a \,\diamond\, a \,=\, a$
\end{Definition}

\begin{Remark}
\cite{Park}
In intuitionistic fuzzy metric space $X$, \, $\mu\,(\,x \,,\,y \,,\,\cdot\,)$ \,
is non-decreasing and \,$\nu\,(\,x \,,\,y \,,\,\cdot\,)$\, is non-increasing
for all \,$x\,,\,y\,\in\,X$.
\end{Remark}

\begin{Definition}
\cite{Veeramani}
 A sequence $\{\,x_{\,n}\,\}_{\,n}$ in a intuitionistic fuzzy metric space
 is said to be a \textit{\textbf{Cauchy sequence}} if and only if for each
 \,$r\,\in\,(\,0\,,\,1\,)$\, and \,$t>0$\, there exists \,$n _{\,0}\,\in\,N$\,
 such that \,$\mu\,(\,x_{n} \,,\,x_{m} \,,\, t\,) \;>\;1\,-\,r $\, and
 \,$\nu\,(\,x_{n} \,,\,x_{m} \,,\, t\,) \;<\;r$\, for all
 $\;n\,,\,m\;\geq\;n_{\,0}.$\, \\A sequence $\{\,x_{\,n}\,\}$ in
 a intuitionistic fuzzy metric space is said to converge to \,$x\,\in\,X$\,
 if and only if for each \,$r\,\in\,(\,0\,,\,1\,)$\, and \,$t>0$\,
 there exists \,$n _{\,0}\,\in\,N$\, such that
 \,$\mu\,(\,x_{n} \,,\,x \,,\, t\,) \;>\;1\,-\,r $\, and
 \,$\nu\,(\,x_{n} \,,\,x \,,\, t\,) \;<\;r$\, for all $\;n\,,\,m\;\geq\;n_{\,0}.$\,
\end{Definition}

\begin{Note}
\cite{Samanta}
A sequence $\{\,x_{\,n}\,\}_{\,n}$ in an intuitionistic fuzzy metric space
is a Cauchy sequence if and only if \[\mathop {\lim }\limits_{n\,\, \to
\,\,\infty } \,\mu\,(\,x_{n} \,,\,x_{n\,+\,p} \,,\, t\,)\;=\;1\;\; and
\;\; \mathop {\lim }\limits_{n\,\, \to \,\,\infty } \,\nu\,(\,x_{n} \,,
\,x_{n\,+\,p} \,,\, t\,)\;=\;0\;\]
A sequence $\{\,x_{\,n}\,\}_{\,n}$ in an intuitionistic fuzzy metric space converges to
\,$x\,\in\,X$\, if and only if \[\mathop {\lim }\limits_{n\,\, \to \,\,\infty }
\,\mu\,(\,x_{n} \,,\,x \,,\, t\,)\;=\;1\;\; and \;\; \mathop {\lim }\limits_{n\,\,
\to \,\,\infty } \,\nu\,(\,x_{n} \,,\,x \,,\, t\,)\;=\;0\;\]
\end{Note}

\begin{Definition}
\cite{Abdul}
 Let $(\,X\,,\,A\,)$ be a intuitionistic fuzzy metric space. We will say
 the mapping \,$f : X \rightarrow X$\, is \textbf{t-uniformly continuous}
 if for each \,$\varepsilon,$\, with \,$0\, < \,\varepsilon \,< \,1,$\,
 there exists  \,$0 \,<\, r \,<\, 1,$\, such that
 \,$\mu\,(\,x \,,\, y \,,\, t\,) \,\geq\,1\,-\,r\, $\, and
 \,$\nu\,(\,x \,,\, y \,,\, t\,) \,\leq\,r\, $ implies
 $\mu\,\left(\,f(x) \,,\, f(y) \,,\, t\,\right) \,\geq\,1\,-\,\varepsilon \, $
 and \,$\nu\,\left(\,f(x) \,,\, f(y) \,,\, t\,\right) \,\leq\,\varepsilon \, $
 for each \,$x , y\,\in\,X$\, and \,$t \,> \,0.$
\end{Definition}

\begin{Definition}
\cite{T.Samanta}
Let \,$(\,X \,,\, A\,)$\, be \textit{\textbf{IFMS}} and $\;T:X\rightarrow X$.
T is said to be TS-IF contractive mapping if there exists \,$k\;\in\;(0 \,,\, 1)$
such that \[k\;\mu\,\left(\,T(x) \;,\; T(y) \;,\; t\,\right)\;\geq\;
\mu\,(\,x \;,\; y \;,\; t\,)\;\] \[ and  \hspace{0.8cm}\frac{1}{k}\;
\nu\,\left(\,T(x) \;,\; T(y) \;,\; t\,\right)\;\leq\;\nu\,(\,x \;,\;
y \;,\; t\,) \;\;\forall \;\;t \,>\, 0.\]
\end{Definition}

\begin{Proposition}\label{P1}
Let \,$(\,X \,,\, A\,)$\, be a intuitionistic fuzzy metric space.
If \,$f:X\rightarrow\;X$ is TS-IF contractive then f is t-uniformly continuous.
\end{Proposition}

{\bf Proof.} Obvious

\begin{Theorem}
\cite{T.Samanta}
Let \,$(\,X \,,\, A\,)$\, be a complete \textit{\textbf{IFMS}} and
$\,T : X \rightarrow X$ be TS-IF contractive mapping with k its
contraction constant. Then T has a unique fixed point.
\end{Theorem}


\section{ Fixed-point theorems}


\begin{Definition}
Let \,$(\,X \,,\, A\,)$\, be an \textit{\textbf{IFMS}} , $ x\,\in\,X,$\,
$r\,\in\,(0 \,,\, 1)$ , \,$t \,>\, 0$,
 \\${\hspace{2.0cm}}B(\,x \,,\, r \,,\, t\,)\,=\,\left\{\,y\,\in \,X\,/\,
 \mu(\,x \,,\, y \,,\, t\,) \,>\, 1\,-\,r \,,\, \nu(\,x \,,\, y \,,\, t\,)
 \,<\, r\,\right\}$.
 \\ Then \,$B(\,x \,,\, r \,,\, t\,)$\, is called an \textbf{open ball}
 centered at x of radius r w.r.t. t.
\end{Definition}

\begin{Definition}
Let \,$(\,X \,,\, A\,)$\, be an \textit{\textbf{IFMS}} and \,$P\,\subseteq \,X$\,.
P is said to be a \textbf{closed set} in \,$(\,X \,,\, A\,)$\, if and only if
any sequence \,$\{\,x_{n}\}$\, in P converges to \,$x\,\in P$\, i.e, iff.
$\mathop {\lim }\limits_{n\; \to \;\infty } \,\mu\,(\,x_{n} \,,\, x \,,\, t\,)\,=\,1$\,
and \,$\mathop {\lim }\limits_{n\; \to \;\infty } \,\nu\,(\,x_{n}
\,,\,x \,,\, t\,)\,=\,0\;$  $\Rightarrow$ $x\in\,P$.
\end{Definition}

\begin{Definition}\label{D1}
Let \,$(\,X \,,\, A\,)$\, be an \textit{\textbf{IFMS,}} $ x\,\in\,X,$\,
$r\,\in\,(0 \,,\, 1)$ , \,$t \,>\, 0$,
 \\${\hspace{2.0cm}}S(\,x \,,\, r \,,\, t\,)\,=\,\left\{\,y\,\in \,X\,/\,
 \mu(\,x \,,\, y \,,\, t\,) \,>\, 1\,-\,r \,,\, \nu(\,x \,,\, y \,,\, t\,)
 \,<\, r\,\right\}$.
 \\Hence \,$S(\,x \,,\, r \,,\, t\,)$\, is said to be a \textbf{closed ball}
 centered at x of radius r w.r.t. t iff. any sequence $\{\,x_{n}\}$ in
 \,$S(\,x \,,\, r \,,\, t\,)$\, converges to y then $y\in\,\,S(\,x \,,\, r \,,\, t\,)$.
\end{Definition}

\begin{Theorem}
$($Contraction on a closed ball$)$ :- Suppose \,$(\,X \,,\, A\,)$\,
is a complete IFMS. Let $T:X\,\rightarrow\,X$ be TS-IF contractive
mapping on \,$S(\,x \,,\, r \,,\, t\,)$ \,with contraction constant k.
Moreover, assume that
\[k\,{\mu\,\left(\,x \,,\,T(x) \,,\, t\,\right)}
\,>\, (1 \,-\, r)\;\; and \;\; \frac{1}{k}
\;{\nu\,\left(\,x \,,\,T(x) \,,\, t\,\right)}\,<\,r\]
Then T has unique  fixed point in \,$S(\,x \,,\, r \,,\, t\,)$.
\end{Theorem}

{\bf Proof.}  Let $x_{1}\,=\,T(x) \,,\, x_{2}\,=\,T(x_{1})\,=\,T^{\,2}(x)
\,,\; \cdots \;,\, x_{n}\,=\,T(x_{n-1})$ \\i.e, \,$\,x_{n}\,=\,T^{\,n}(x)$\,
for all $n\,\in\,N.\;$ Now \[k\,{\mu\,\left(\,x \,,\,T(x) \,,\, t\,\right)}
\,>\, (1 \,-\, r){\hspace{2.6cm}}\]
\[\Rightarrow\;\;{\mu\,\left(\,x \,,\,T(x) \,,\, t\,\right)} \,>\,
\frac{(1 \,-\, r)}{k} \,>\, (1 \,-\, r)\hspace{1.5cm} \]
\[\Rightarrow\;\;{\mu\,\left(\,x \,,\,x_1 \,,\, t\,\right)} \,>\,
(1 \,-\, r) \hspace{1.5cm} \cdots \hspace{1.5cm}(\,i\,)\]
Again,
\[\frac{1}{k}\; {\nu\,\left(\,x \,,\,T(x) \,,\, t\,\right)}\,<\,
r {\hspace{2.5cm}}\]
\[\Rightarrow\;\;{\nu\,\left(\,x \,,\,T(x) \,,\, t\,\right)} \,<\,
 r\,k\,<\,r\hspace{3.2cm} \]
\[\Rightarrow\;\;{\nu\,(\,x \,,\, x_1 \,,\, t\,)}\,<\,r  \hspace{1.0cm}
\cdots  \hspace{2.5cm} (\,ii\,) {\hspace{0.1cm}}\]
$(\,i\,)$ and $(\,ii\,)$ $\;\Rightarrow\; x_{1}\,\in\,S(\,x \,,\, r \,,\, t\,).$
\\\\Assume that \,$x_{1} \,,\, x_{2} \,,\; \cdots \;,\, x_{n \,-\, 1}
\,\in\,S(\,x \,,\, r \,,\, t\,)$.
We show that $x_{n}\,\in\,S(\,x \,,\, r \,,\, t\,).$
\[k\;\mu\,(\,x_{1} \;,\; x_{2}
\;,\; t\,)\;=\;k\;\mu\,\left(\,T(x) \;,\; T(x_{1}) \;,\; t\,\right) {\hspace{3.5cm}}\]
\[\;\geq\;\mu\,(\,x \;,\; x_{1}
\;,\; t\,) {\hspace{2.2cm}}\]
\[\Rightarrow\;\mu\,(\,x_{1} \;,\; x_{2}
\;,\; t\,)\,>\,\frac{(1-r)}{k}\,>\,(1-r)\hspace{4.8cm}\]
\[k\;\mu\,(\,x_{2} \;,\; x_{3}
\;,\; t\,)\;=\;k\;\mu\,\left(\,T(x_{1}) \;,\; T(x_2) \;,\; t\,\right) {\hspace{3.5cm}}\]
\[\;\geq\;\mu\,(\,x_{1} \;,\; x_{2}
\;,\; t\,) {\hspace{2.2cm}}\]
\[\Rightarrow\;\mu\,(\,x_{2} \;,\; x_{3}
\;,\; t\,)\,\geq\,\frac{1}{k}\,\mu\,(\,x_{1} \;,\; x_{2}
\;,\; t\,)\hspace{5.2cm}\]
\[>\;\,\frac{1-r}{k}\,>\,(1-r)\hspace{1.7cm}\]
Again,
\[\frac{1}{k}\;\nu\,(\,x_{1} \;,\; x_{2}
\;,\; t\,)\;=\;\frac{1}{k}\;\nu\,\left(\,T(x) \;,\; T(x_{1}) \;,\; t\,\right)
 {\hspace{3.5cm}}\]
\[\;\leq\;\nu\,(\,x \;,\; x_{1}
\;,\; t\,) {\hspace{2.2cm}}\]
\[\Rightarrow\;\nu\,(\,x_{1} \;,\; x_{2}
\;,\; t\,)\,\leq\,k\,\nu\,(\,x \;,\; x_{1}
\;,\; t\,)\,<\,k\,r\,<\,r{\hspace{3.2cm}}\]
\[\frac{1}{k}\;\nu\,(\,x_{2} \;,\; x_{3}
\;,\; t\,)\;=\;\frac{1}{k}\;\nu\,\left(\,T(x_{1}) \;,\; T(x_{2}) \;,\;
 t\,\right) {\hspace{3.5cm}}\]
\[\;\leq\;\nu\,(\,x_{1} \;,\; x_{2}
\;,\; t\,) {\hspace{2.2cm}}\]
\[\Rightarrow\;\nu\,(\,x_{2} \;,\; x_{3}
\;,\; t\,)\,\leq\,k\,\nu\,(\,x_{1} \;,\; x_{2}
\;,\; t\,)\,<\,k\,r\,<\,r{\hspace{3.2cm}}\]
Similarly it can be shown that ,\\
$\mu\,(\,x_{3} \,,\, x_{4} \,,\, t\,) \;>\; 1 \,-\, r \;,\;
\nu\,(\,x_{3} \,,\, x_{4} \,,\, t\,) \;<\; r$ \,,\,
$\; \cdots\; ,\, \mu\,(\,x_{n \,-\, 1} \,,\, x_{n} \,,\, t\,) \;>\; 1 \,-\, r$ and
$\;\nu\,(\,x_{n \,-\, 1} \,,\, x_{n} \,,\, t\,) \;<\; r.$ \\Thus, we see that ,
\[\mu\,(\,x \,,\, x_{n} \,,\, t\,)\;\geq\;\mu\,\left(\,x \,,\, x_{1}
\,,\, \frac{t}{n}\,\right)\;\ast\;\mu\,\left(\,x_{1} \,,\, x_{2} \,,\,
\frac{t}{n}\,\right)\;\ast\;\cdots\;\ast\;\mu\,\left(x_{n \,-\, 1}
\,,\, x_{n} \,,\, \frac{t}{n}\,\right)\]
\[>\; (1 \,-\, r)\;\ast\;(1 \,-\, r)\;\ast\;\cdots\;\ast\;(1 \,-\, r)
\;=\; 1 \,-\, r\;\;\]
\[i.e. \,, \hspace{0.5cm}\mu\,(\,x \,,\, x_{n} \,,\, t\,) \;>\; 1
\,-\, r {\hspace{9.5cm}}\]
\[\nu\,(\,x \,,\, x_{n} \,,\, t\,)\;\leq\;\nu\,\left(\,x \,,\, x_{1}
\,,\, \frac{t}{n}\,\right)\;\diamond\;\nu\,\left(\,x_{1} \,,\, x_{2} \,,\,
\frac{t}{n}\,\right)\;\diamond\;\cdots\;\diamond\;\nu\,\left(x_{n \,-\, 1}
\,,\, x_{n} \,,\, \frac{t}{n}\,\right)\]
\[\;<\;r\;\diamond\;r\;\diamond\;\cdots\;\diamond\;r\;=\;r {\hspace{4.6cm}}\]
Thus , $\;\mu\,(\,x \,,\, x_{n} \,,\, t\,) \;>\; 1 \,-\, r \;and\;
\nu\,(\,x \,,\, x_{n} \,,\, t\,) \;<\; r$
\\$\Rightarrow\;x_{n}\;\in\;S(\,x \,,\, r \,,\, t\,)$
\\Hence, by the theorem 3.10\cite{T.Samanta} and the definition\ref{D1},
 T has unique fixed point in \,$S(\,x \,,\, r \,,\, t\,)$.

\begin{Note}
It follows from the proof of Theorem 3.10\cite{T.Samanta} that for
any $\;x\in\;X$ the sequence of iterates \,$\{\,T^{\,n}(x)\}$\,
converges to the fixed point of T.
\end{Note}

\begin{Lemma}\label{L2}
Let \,$(\,X \,,\, A\,)$\, be \textit{\textbf{IFMS}} and $T:\,X\rightarrow\,X$
be t-uniformly  continuous on X. If $\;x_{n}\;\rightarrow\;x\; $ as
$\; n\;\rightarrow\infty\;$\, in \,$(\,X \,,\, A\,)$\, then
$\;T(x_{n})\;\rightarrow\;T(x)\;$ as $\; n\;\rightarrow\infty\;$
in \,$(\,X \,,\, A\,)$ .
\end{Lemma}
{\bf Proof.} \, Proof directly follows from the definitions of
t-uniformly continuity and convergence of a sequence in a \textit{\textbf{IFMS}}.

\begin{Lemma}\label{L1}
Let \,$(\,X \,,\, A\,)$\, be \textit{\textbf{IFMS}}.
If $\;x_{n}\;\rightarrow\;x$ and $\;y_{n}\;\rightarrow\;y$ in
\,$(\,X \,,\, A\,)$\, then $\mu\,(\,x_{n} \,,\, y_{n} \,,\, t\,)
\;\rightarrow\;\mu\,(\,x \,,\, y \,,\, t\,) \;and\; \nu\,(\,x_{n}
\,,\, y_{n} \,,\, t\,)\;\rightarrow\;\nu\,(\,x \,,\, y \,,\, t\,)
\;as \;n\;\rightarrow\infty \;\; for\; all\;\; t \,>\, 0\; \;in \;\;R\;.$
\end{Lemma}

{\bf Proof.} \,\,We have,
\[ {\hspace{1.0cm}}\mathop {\lim }\limits_{n\;\, \to \;\,\infty }
\;\mu\,(\,x_{n} \;,\,x \;,\; t\,)\;=\;1\; , \;\mathop {\lim }\limits_{n\;\,
 \to \;\,\infty } \;\nu\,(\,x_{n} \;,\,x \;,\; t\,)\;=\;0\;\]
\[ and \;\;\mathop {\lim }\limits_{n\;\, \to \;\,\infty } \;
\mu\,(\,y_{n} \;,\,y \;,\; t\,)\;=\;1\; , \;\mathop {\lim }\limits_{n\;\,
\to \;\,\infty } \;\nu\,(\,y_{n} \;,\,y \;,\; t\,)\;=\;0\;\]
\[ \mu\,(\,x_{n} \,,\, y_{n} \,,\, t\,)\;\geq\;\mu\,\left(\,x_{n} \,,
\, x \,,\, \frac{t}{2}\,\right)\;\ast\;\mu\,\left(\,x \,,\,y_{n} \,,\,
\frac{t}{2}\,\right)\]
\[{\hspace{4.8cm}}\;\geq\;\mu\,\left(\,x_{n} \,,\, x \,,\, \frac{t}{2}\,
\right)\;\ast\;\mu\,\left(\,x \,,\, y \,,\, \frac{t}{4}\,\right)\;\ast\;\mu\,
\left(\,y \,,\, y_{n} \,,\, \frac{t}{4}\,\right)\]
\\${\hspace{0.5cm}}\Rightarrow\;\mathop {\lim }\limits_{n\;\, \to \;\,\infty }
\;\mu\,(\,x_{n} \;,\,y_{n} \;,\; t\,)\;\geq\;\mu\,(\,x\;,\;y\;,\;t\,)\;$
\[ \mu\,(\,x \,,\, y \,,\, t\,)\;\geq\;\mu\,\left(\,x \,,\, x_{n} \,,\,
\frac{t}{2}\,\right)\;\ast\;\mu\,\left(\,x_{n} \,,\, y \,,\, \frac{t}{2}\,\right)
{\hspace{2.5cm}}\]
\[ {\hspace{2.7cm}}\;\geq\;\mu\,\left(\,x \,,\, x_{n} \,,\, \frac{t}{2}\,\right)\;
\ast\;\mu\,\left(\,x_{n} \,,\, y_{n} \,,\, \frac{t}{4}\,\right)\;\ast\;\mu\,
\left(\,y_{n} \,,\, y \,,\, \frac{t}{4}\,\right)\]
${\hspace{0.5cm}}\Rightarrow\;\mu\,(\,x \;,\,y \;,\; t\,)\;\geq\;\;
\mathop {\lim }\limits_{n\;\, \to \;\,\infty } \;\mu\,(\,x_{n}\;,\;y_{n}\;,\;t\,)\;$
$\;\forall\;t\,\,>\,0.$
\\Then,\\ ${\hspace{0.9cm}}\;\mathop {\lim }\limits_{n\;\, \to \;\,\infty }
\;\mu\,(\,x_{n} \;,\;y_{n} \;,\; t\,)\;=\;\mu\,(\,x\;,\;y\;,\;t\,)\;$
for all $\;t\;>\,0$,
\\\\Similarly, $\;\mathop {\lim }\limits_{n\;\, \to \;\,\infty } \;
\nu\,(\,x_{n} \;,\;y_{n} \;,\; t\,)\;=\;\nu\,(\,x\;,\;y\;,\;t\,)\;$
for all $\;t\;>\,0$.

\begin{Theorem}
Let \,$(\,X \,,\, A\,)$\, be a complete \textit{\textbf{IFMS}} and
 \,$T:X\rightarrow X$\, be a t-uniformly continuous on X. If for same
 positive integer m, \,$T^{\,m}\,$ is a TS-IF contractive mapping with
 k its contractive constant then $T$ has a unique fixed point in \,$X$.
\end{Theorem}

{\bf Proof.} \,Let $B\;=\;T^{\,m}$,\,  $n$ \,be an arbitrary but fixed
positive integer and \,$x\;\in\;X$. \\we now show that $\;B^{\,n}\,T(x)
\;\rightarrow\; B^{\,n}\,(x)$ in \,$(\,X \,,\, A\,)$\,.
\\Now,
\[k\;\mu\,\left(\,B^{\,n}T(x) \;,\; B^{\,n}(x) \;,\; t\,\right)\;=\;
k\;\mu\,\left(\,B(B^{\,n \,-\, 1}\,T(x)) \;,\; B(B^{\,n \,-\, 1}(x))
\;,\; t\,\right)  {\hspace{3.5cm}}\]
\[\;{\hspace{1.5cm}}\geq\;{\mu\,(\,B^{n \,-\, 1}\,T(x) \;,\; B^{n \,-
\, 1}(x) \;,\; t\,)}\]
\[i.e,\;\mu\,\left(\,B^{\,n}T(x) \;,\; B^{\,n}(x) \;,\; t\,\right)\;\geq\;
\,\frac{1}{k}\,\;{\mu\,(\,B^{n \,-\, 1}\,T(x) \;,\; B^{n \,-\, 1}(x) \;,\; t\,)}
 {\hspace{2.8cm}} \]
\[{\hspace{4.5cm}}\;=\;\frac{1}{k}\,\;\mu\,\left(\,B(B^{\,n \,-\, 2}\,T(x))
 \;,\; B(B^{\,n \,-\, 2}(x)) \;,\; t\,\right) \]
\[{\hspace{5.2cm}}\;\geq\;\frac{1}{k^{2}}\,\;{\mu\,(\,B^{n \,-\, 2}\,T(x)
\;,\; B^{n \,-\, 2}(x) \;,\; t\,)} {\hspace{2.2cm}}\]
\[\;\geq\;{\hspace{1.0cm}}\cdots{\hspace{1.1cm}}\]
\[{\hspace{0.6cm}}\;\geq\;\frac{1}{k^{n}}\,\;\mu\,(\,T(x) \;,\; x \;,\; t\,)\]
\[ \Rightarrow\;\;\mathop {\lim }\limits_{n\; \to \;\infty}
\mu\,\left(\,B^{\,n}T(x) \;,\; B^{\,n}(x) \;,\; t\,\right)\;
\geq\;\mathop {\lim }\limits_{n\; \to \;\infty }\;\frac{1}{k^{n}}\;
\mu\,(\,T(x) \;,\; x \;,\; t\,)\]
\[\Rightarrow \;\;\mathop {\lim }\limits_{n\; \to \;\infty} \mu\,
\left(\,B^{\,n}T(x) \;,\; B^{\,n}(x) \;,\; t\,\right)\;=\;1 {\hspace{4.4cm}}\]
\\Similarly,
${\hspace{0.5cm}}\mathop {\lim }\limits_{n\;\, \to \;\,\infty }
\nu\,\left(\,B^{\,n}\,T(x) \;,\; B^{\,n}(x) \;,\; t\,\right)\;=\;0,\;$
for all $\; t\;>\;0$. \\
Thus, $\;B^{\,n}\,T(x) \;\rightarrow\; B^{\,n}\,(x)$ in \,$(\,X \,,\, A\,)$\,.
\\Again, by the theorem 3.10\cite{T.Samanta}, we see that B has a unique
fixed point \,$y$(say), and from the note [3.5], it follows that
$\;B^{\,n}(x)\,\rightarrow \,y$ as $\;n\;\rightarrow\;\infty\;$
in \,$(\,X \,,\, A\,)$\,.
\\\\Since T is t-uniformly continuous on X, it follows from the above
lemma[3.6] that $\,B^{\,n}\,T(x)\;=\;T\;B^{\,n}(x)\;\rightarrow\;T(y)$
as $n\;\rightarrow\;\infty$ in \,$(\,X \,,\, A\,)$\,.
\\\\ Again since $\mathop {\lim }\limits_{n\;\, \to \;\,\infty }\mu\,
\left(\,B^{\,n}\,T(x) \;,\; B^{\,n}(x) \;,\; t\,\right)\;=\;1\;$ and
$\;\mathop {\lim }\limits_{n\;\, \to \;\,\infty }\nu\,\left(\,B^{\,n}\,T(x)
\;,\; B^{\,n}(x) \;,\; t\,\right)\\=\;0$, \,we have  by the lemma [\ref{L1}]
\\$\mathop {\lim }\limits_{n\;\, \to \;\,\infty }\mu\,\left(\,T(y) \;,
\; y \;,\; t\,\right)\;=\;1$ and $\mathop {\lim }\limits_{n\;\, \to
\;\,\infty }\nu\,\left(\,T(y) \;,\; y \;,\; t\,\right)\;=\;0,\;$ for
all $\;t\;>\;0$,
\\$i.e. \,, \;\;\mu\,\left(\,T(y) \;,\; y \;,\; t \,\right)\;=\;1\;$
and $\;\nu\,\left(\,T(y) \;,\; y \;,\; t\,\right)\;=\;0,\;$ for all $\;t\;>\;0$.
\\$\Rightarrow\;T(y)\;=\;y\;$
$\;\Rightarrow\;y$ is a fixed point of T. \\If $y^{'}$ is a fixed point of T,
$\,i.e.\,, \;T(y^{'})\;=\;y^{'}$, then $\;T^{\,m}(y^{'})\;=\;T^{\,m \,-
\, 1}(T(y^{'}))\;=\;T^{\,m \,-\, 1}(y^{'})\;=\; \cdots \;=\;y^{'}$
$\Rightarrow\;B(y^{'})\;=\;y^{'}\;\Rightarrow\;y^{'}\,$ is a fixed point of B.
 \\But $y$ is the unique fixed point of B, therefore $y\;=\;y^{'}$ which
 implies that $y$ is the unique fixed point of T. This completes the proof.

\begin{Definition}
Let \,$(\,X \,,\, A\,)$\, be a \textit{\textbf{IFMS}} and \,$T:X\rightarrow X$.
For \,$\varepsilon\;>\;0$\, and \,$0 \;<\; \lambda \;<\; 1$\, , $T$\,
is said to be \,$(\,\varepsilon \,,\, \lambda\,)$\,
\textit{\textbf{IF-uniformly locally contractive}}
if \[\mu\,(\,x\;,\;y\;,\;t\,) \;>\; \varepsilon \;\Longrightarrow\;
\lambda\;\mu\,(\,Tx\;,\;Ty\;,\;t\,) \;>\; \mu\,(\,x\;,\;y\;,\;t\,)\]
\[\nu\,(\,x\;,\;y\;,\;t\,) \;<\; 1 \;-\; \varepsilon \;\Longrightarrow\;
\frac{1}{\lambda}\;\nu\,(\,Tx\;,\;Ty\;,\;t\,) \;<\; \nu\,(\,x\;,\;y\;,\;t\,)
\hspace{0.9cm}\]
\end{Definition}

\begin{Definition}
Let \,$0 \;<\; \eta \;<\; 1$\, and \,$(\,X \,,\, A\,)$\, be a \textit{\textbf{IFMS}}.
Then \,$(\,X \,,\, A\,)$\, is said to be \textit{\textbf{IF}}\,
\,$\eta\,-$\,\textit{\textbf{chainable}} space if for every
\,$a \,,\, b \; \in \;X$\, there exist a finite set of points
\,$a \,=\, x_{\,0} \;,\; x_{\,1} \;,\; \cdots \;,\; x_{\,n} \,=\, b$\, such that
\[\mu\,(\,x_{\,i \,-\, 1}\;,\;x_{\,i}\;,\;t\,) \;>\; \eta \hspace{0.3cm}
and \hspace{0.3cm} \nu\,(\,x_{\,i \,-\, 1}\;,\;x_{\,i}\;,\;t\,) \;<\; 1
\;-\; \eta \; \; ,\hspace{0.3cm} i \,=\, 1 \,,\, 2 \,,\; \cdots \;,\, n\]
\end{Definition}

\begin{Theorem}
Let \,$(\,X \,,\, A\,)$\, be a complete \textit{\textbf{IFMS}} and
IF\,$\varepsilon\,-$\,chainable space. If \,$T:X\rightarrow X$\, is
\,$(\,\varepsilon \,,\, \lambda\,)$\, IF-uniformly locally contractive
then $T$ has a fixed point in \,$X$.
\end{Theorem}

{\bf Proof.} \,Let \,$x$\, be an arbitrary but fixed point of \,$X$.
If \,$T x \,=\, x$\, then a fixed point is assured. We assume therefore
that \,$T x \;\neq\; x$. Since \,$X$\, is IF\,$\varepsilon\,-$\,chainable
space, there exists a finite set of points \,$x \,=\, x_{\,0} \;,\; x_{\,1}
\;,\; \cdots \;,\; x_{\,n} \,=\, T x$\, such that

\[\mu\,(\,x_{\,i \,-\, 1}\;,\;x_{\,i}\;,\;t\,) \;>\; \varepsilon \hspace{0.3cm}
and \hspace{0.3cm} \nu\,(\,x_{\,i \,-\, 1}\;,\;x_{\,i}\;,\;t\,) \;<\; 1 \;-\;
\varepsilon \; \; ,\hspace{0.3cm} i \,=\, 1 \,,\, 2 \,,\; \cdots \;,\, n\]

Again, since \,$T$\, is \,$(\,\varepsilon \,,\, \lambda\,)$\, IF-uniformly
locally contractive, we have
\[ \mu\,(\,x_{\,i \,-\, 1}\;,\;x_{\,i}\;,\;t\,) \;>\; \varepsilon \;
\Longrightarrow \; \lambda\;\mu\,(\,T x_{\,i \,-\, 1}\;,\;T x_{\,i}\;,\;t\,)
 \;>\; \mu\,(\,x_{\,i \,-\, 1}\;,\;x_{\,i}\;,\;t\,) \;>\; \varepsilon\]
\[i.e. , \hspace{0.5cm} \mu\,(\,T x_{\,i \,-\, 1}\;,\;T x_{\,i}\;,\;t\,)
\;>\; \frac{\varepsilon}{\lambda} \;>\; \varepsilon \;;\hspace{6.5cm}\]

\[ \nu\,(\,x_{\,i \,-\, 1}\;,\;x_{\,i}\;,\;t\,) \;<\; 1 \;-\; \varepsilon \;
 \Longrightarrow \; \frac{1}{\lambda}\;\nu\,(\,T x_{\,i \,-\, 1}\;,
 \;T x_{\,i}\;,\;t\,) \;<\; \nu\,(\,x_{\,i \,-\, 1}\;,\;x_{\,i}\;,\;t\,)
 \;<\; 1 \;-\; \varepsilon\]
\[i.e. , \hspace{0.5cm} \nu\,(\,T x_{\,i \,-\, 1}\;,\;T x_{\,i}\;,\;t\,)
\;<\; \lambda\,(\,1 \;-\; \varepsilon\,) \;<\; 1 \;-\; \varepsilon \hspace{6.5cm}\]

and therefore, \[ \lambda^{\,2}\;\mu\,(\,T^{\,2} x_{\,i \,-\, 1}\;,
\;T^{\,2} x_{\,i}\;,\;t\,) \;=\; \lambda\;\left(\, \lambda\;
\mu\,(\,T\,(\,T x_{\,i \,-\, 1}\,)\;,\; T\,(\,T x_{\,i}\,)\;,\;t\,) \,\right)\]
\[\hspace{3.7cm} >\; \lambda\;\mu\,(\,T x_{\,i \,-\, 1}\;,\;T x_{\,i}\;,\;t\,)
\;>\; \lambda\;\varepsilon\]
\[\Longrightarrow\;\; \mu\,(\,T^{\,2} x_{\,i \,-\, 1}\;,\;T^{\,2} x_{\,i}\;,\;t\,)
\;>\; \varepsilon \;;\hspace{6.5cm}\]

\[ \frac{1}{\lambda^{\,2}}\;\nu\,(\,T^{\,2} x_{\,i \,-\, 1}\;,
\;T^{\,2} x_{\,i}\;,\;t\,) \;=\; \frac{1}{\lambda}\;\left(\, \frac{1}{\lambda}\;
\nu\,(\,T\,(\,T x_{\,i \,-\, 1}\,)\;,\; T\,(\,T x_{\,i}\,)\;,\;t\,) \,\right)\]
\[\hspace{5.0cm} <\; \frac{1}{\lambda}\;\nu\,(\,T x_{\,i \,-\, 1}\;,\;T x_{\,i}\;,
\;t\,) \;<\; \frac{1}{\lambda}\;(\,1 \;-\; \varepsilon\,)\]
\[{\hspace{0.2cm}}\Longrightarrow\;\; \nu\,(\,T^{\,2} x_{\,i \,-\, 1}\;,
\;T^{\,2} x_{\,i}\;,\;t\,) \;<\; 1 \;-\; \varepsilon \hspace{6.5cm}\]

In the similar way we have,
\[ \lambda^{\,3}\;\mu\,(\,T^{\,3} x_{\,i \,-\, 1}\;,\;T^{\,3} x_{\,i}\;,
\;t\,) \;=\; \lambda^{\,2}\;\left(\, \lambda\;\mu\,(\,T\,(\,T^{\,2}
x_{\,i \,-\, 1}\,)\;,\; T\,(\,T^{\,2} x_{\,i}\,)\;,\;t\,) \,\right)\]
\[\hspace{3.7cm} >\; \lambda^{\,2}\;\mu\,(\,T^{\,2} x_{\,i \,-\, 1}\;,
\;T^{\,2} x_{\,i}\;,\;t\,) \;>\; \lambda^{\,2}\;\varepsilon\]
\[\Longrightarrow\;\; \mu\,(\,T^{\,3} x_{\,i \,-\, 1}\;,\;T^{\,3}
x_{\,i}\;,\;t\,) \;>\; \varepsilon \hspace{7.5cm}\]
\[\cdots\]
\[ \lambda^{\,m}\;\mu\,(\,T^{\,m} x_{\,i \,-\, 1}\;,\;T^{\,m} x_{\,i}\;,\;t\,)
\;=\; \lambda^{\,m \,-\, 1}\;\left(\, \lambda\;\mu\,(\,T\,(\,T^{\,m \,-\, 1}
x_{\,i \,-\, 1}\,)\;,\; T\,(\,T^{\,m \,-\, 1} x_{\,i}\,)\;,\;t\,) \,\right)\]
\[\hspace{5.2cm} >\; \lambda^{\,m \,-\, 1}\;\mu\,(\,T^{\,m \,-\, 1}
x_{\,i \,-\, 1}\;,\;T^{\,m \,-\, 1} x_{\,i}\;,\;t\,) \;>\; \lambda^{\,m \,-\, 1}
\;\varepsilon\]
\[\Longrightarrow\;\; \mu\,(\,T^{\,m} x_{\,i \,-\, 1}\;,\;T^{\,m} x_{\,i}\;,\;t\,)
\;>\; \varepsilon \;;\hspace{7.5cm}\]

\[ \frac{1}{\lambda^{\,3}}\;\nu\,(\,T^{\,3} x_{\,i \,-\, 1}\;,\;T^{\,3}
x_{\,i}\;,\;t\,) \;=\; \frac{1}{\lambda^{\,2}}\;\left(\, \frac{1}{\lambda}
\;\nu\,(\,T\,(\,T^{\,2} x_{\,i \,-\, 1}\,)\;,\; T\,(\,T^{\,2} x_{\,i}\,)\;,
\;t\,) \,\right)\]
\[\hspace{5.2cm} <\; \frac{1}{\lambda^{\,2}}\;\nu\,(\,T^{\,2} x_{\,i \,-\, 1}
\;,\;T^{\,2} x_{\,i}\;,\;t\,) \;<\; \frac{1}{\lambda^{\,2}}\;(\,1 \;-\;
\varepsilon\,)\]
\[\Longrightarrow\;\; \nu\,(\,T^{\,3} x_{\,i \,-\, 1}\;,\;T^{\,3} x_{\,i}\;,
\;t\,) \;<\; (\,1 \;-\; \varepsilon\,) \hspace{7.5cm}\]
\[\cdots\]
\[ \frac{1}{\lambda^{\,m}}\;\nu\,(\,T^{\,m} x_{\,i \,-\, 1}\;,\;T^{\,m}
 x_{\,i}\;,\;t\,) \;=\; \frac{1}{\lambda^{\,m \,-\, 1}}\;
 \left(\, \frac{1}{\lambda}\;\nu\,(\,T\,(\,T^{\,m \,-\, 1} x_{\,i \,-\, 1}\,)
 \;,\; T\,(\,T^{\,m \,-\, 1} x_{\,i}\,)\;,\;t\,) \,\right)\]
\[\hspace{3.7cm} <\; \frac{1}{\lambda^{\,m \,-\, 1}}\;\nu\,(\,T^{\,m \,-\, 1}
 x_{\,i \,-\, 1}\;,\;T^{\,m \,-\, 1} x_{\,i}\;,\;t\,) \;<\;
 \frac{1}{\lambda^{\,m \,-\, 1}}\;(\,1 \;-\; \varepsilon\,)\]
\[\Longrightarrow\;\; \nu\,(\,T^{\,m} x_{\,i \,-\, 1}\;,\;T^{\,m} x_{\,i}\;,
\;t\,) \;<\; 1 \;-\; \varepsilon \hspace{7.5cm}\]
Now,\[ \mu\,(\,T^{\,m}\,x \;,\;T^{\,m \,+\, 1}\,x  \;,\; t\,) \;=\;
\mu\,(\,T^{\,m}\,x_{\,0} \;,\;T^{\,m}\,x_{\,n}  \;,\; t\,) \hspace{3.5cm}\]
\[ \hspace{5.1cm} \geq \; \;\left(\, \mu\,\left(\,T^{\,m}\,x_{\,0} \;,
\;T^{\,m}\,x_{\,1}  \;,\; \frac{t}{n}\,\right) \;\ast\;
\mu\,\left(\,T^{\,m}\,x_{\,1} \;,\;T^{\,m}\,x_{\,2}  \;,\;
\frac{t}{n}\,\right)\,\right. \]
\[\left. \hspace{6.9cm} \ast \hspace{0.5cm} \cdots \hspace{0.5cm} \ast \;
\mu\,\left(\,T^{\,m}\,x_{\,n \,-\, 1} \;,\;T^{\,m}\,x_{\,n}  \;,\;
\frac{t}{n}\,\right)\,\right)\] \[ >\; \varepsilon \hspace{3.0cm}\]
\[ i.e. \hspace{0.5cm} \mu\,(\,T^{\,m}\,x \;,\;T^{\,m \,+\, 1}\,x  \;,\; t\,)
\;>\; \varepsilon \hspace{0.5cm} for \;\; all \;\; t \;>\; 0 \; \; and \;\; for
\;\; all \;\; m\;\in\;\mathbf{N}\;;\hspace{5.5cm} \]
and,\[ \nu\,(\,T^{\,m}\,x \;,\;T^{\,m \,+\, 1}\,x  \;,\; t\,) \;=\;
\nu\,(\,T^{\,m}\,x_{\,0} \;,\;T^{\,m}\,x_{\,n}  \;,\; t\,) \hspace{3.5cm}\]
\[ \hspace{5.1cm} \leq \; \;\left(\, \nu\,\left(\,T^{\,m}\,x_{\,0}
\;,\;T^{\,m}\,x_{\,1}  \;,\; \frac{t}{n}\,\right) \;\diamond\;
\nu\,\left(\,T^{\,m}\,x_{\,1} \;,\;T^{\,m}\,x_{\,2}  \;,\;
\frac{t}{n}\,\right)\,\right. \]
\[\left. \hspace{6.9cm} \diamond \hspace{0.5cm} \cdots \hspace{0.5cm}
\diamond \; \nu\,\left(\,T^{\,m}\,x_{\,n \,-\, 1} \;,\;T^{\,m}\,x_{\,n}
\;,\; \frac{t}{n}\,\right)\,\right)\] \[ <\; 1 \;-\; \varepsilon \hspace{2.0cm}\]
\[ i.e. \hspace{0.5cm} \nu\,(\,T^{\,m}\,x \;,\;T^{\,m \,+\, 1}\,x
\;,\; t\,) \;<\; 1 \;-\; \varepsilon \hspace{0.5cm} for \;\; all \;\;
 t \;>\; 0 \; \; and \;\; for \;\; all \;\; m\;\in\;\mathbf{N}.\hspace{5.5cm} \]

Now, for all \,$t \;>\; 0$\, and \,$j \;<\; k$\, we have,
\[ \mu\,(\,T^{\,j}\,x \;,\;T^{\,k}\,x  \;,\; t\,) \;\geq\;
\mu\,\left(\,T^{\,j}\,x \;,\;T^{\,j \,+\, 1}\,x  \;,\; \frac{t}{k \,-\, j}\,\right)
\;\;\ast\;\; \mu\,\left(\,T^{\,j \,+\, 1}\,x \;,\;T^{\,j \,+\, 2}\,x
\;,\; \frac{t}{k \,-\, j}\,\right) \]
\[\hspace{4.5cm} \ast \hspace{0.5cm} \cdots \hspace{0.5cm} \ast \;\;
 \mu\,\left(\,T^{\,k \,-\, 1}\,x \;,\;T^{\,k}\,x  \;,\; \frac{t}{k \,-\, j}\,\right)
 \;\;>\;\;\varepsilon \;;\]
\[ \nu\,(\,T^{\,j}\,x \;,\;T^{\,k}\,x  \;,\; t\,) \;\leq\;
\nu\,\left(\,T^{\,j}\,x \;,\;T^{\,j \,+\, 1}\,x  \;,\; \frac{t}{k \,-\, j}\,\right)
\;\;\diamond\;\; \nu\,\left(\,T^{\,j \,+\, 1}\,x \;,\;T^{\,j \,+\, 2}\,x
\;,\; \frac{t}{k \,-\, j}\,\right) \]
\[\hspace{4.5cm} \diamond \hspace{0.5cm} \cdots \hspace{0.5cm} \diamond
\;\; \nu\,\left(\,T^{\,k \,-\, 1}\,x \;,\;T^{\,k}\,x  \;,\;
\frac{t}{k \,-\, j}\,\right) \;\;<\;\;1 \;-\; \varepsilon.\]
$\Longrightarrow \;\; \left\{\,T^{\,j}\,x\,\right\}$\, is a Cauchy
sequence in \,$(\,X \,,\, A\,)$. Since \,$(\,X \,,\, A\,)$\, is
complete, there exists \,$\xi \;\in\; X$\, such that \,$T^{\,i}\,x
\;\longrightarrow\; \xi$\, as \,$i \,\longrightarrow\, \infty$\, in
\,$(\,X \,,\, A\,)$. Again, since \,$T$\, is IF-uniformly locally
contractive, it follows that \,$T$\, is t-uniformly continuous on X
and hence by the lemma(\ref{L2}), we get \[ T\,\xi \;=\; \mathop
{\lim }\limits_{i\;\, \to \;\,\infty }T\,T^{\,i}x \;=\; \mathop
{\lim }\limits_{i\;\, \to \;\, \infty }T^{\,i \;+\; 1}x \;=\; \xi\]
which shows that \,$\xi$\, is a fixed point of \,$T$.

\end{document}